\documentclass[11pt]{article}
\usepackage{mathrsfs}
\usepackage{CJK}
\usepackage{amsmath}
\usepackage{fancyhdr}

\begin{document}
\title{ Some Results On The Jacobian Conjecture And Polynomial Automorphisms}
\author{{ Dan Yan\footnote{E-mail address:yan-dan-hi@163.com}}     \\
\small School of Mathematical Sciences, Graduate University of
Chinese Academy of Sciences,\\
\small Beijing 100049, China\\}
\date{}
\maketitle \begin {abstract} In this paper, we will first show that, the homogeneous polynomials which satisfy the Jacobian condition are injective on the lines that pass through the origin. Secondly, we will show that $F$ and $G'$ are paired, where $F$ is a Druzkowski map and $G'$ is a cubic homogeneous polynomial which related to $F$. Finally, we will find a more exactly bound for the degree of $F^{-1}$, where $F$ is a invertible map.
\end{abstract} {\bf Keywords.} Jacobian Conjecture,
Polynomial mapping, Druzkowski mapping\\
{\bf MSC(2010).} Primary 14E05  Secondary 14A05;14R15 \vskip 2.5mm

Let
$F=(F_1(x_1,\cdots,x_n),\cdots,F_n(x_1,\cdots,x_n))^t:\bf{C}^n\rightarrow\bf{C}^n$
be a polynomial mapping, that is,
$F_i(x_1,\cdots,x_n)\in\bf{C}[x_1,\cdots,x_n]$ for all $1\leq i\leq
n$. Let $JF=(\frac{\partial F_i}{\partial x_j})_{n\times n}$ be the
Jacobian matrix of $F$. The well-known Jacobian Conjecture(JC)
raised by O.H. Keller in 1939 (\cite{1}) states that a polynomial
mapping $F:\bf{C}^n\rightarrow\bf{C}^n$ is invertible if the
Jacobian determinant $|JF|$ is a nonzero constant. This conjecture
has being attacked by many people from various research fields and
remains open even when $n=2$! (Of course, a positive answer is
obvious when $n=1$ by the elements of linear algebra.)  See \cite{2}
and \cite{3} and the references therein for a wonderful 70-years
history of this famous conjecture. It can be easily seen that JC
is true if JC holds for all polynomial mappings whose Jacobian
determinant is 1. We make use of this convention in the present
paper. Among the vast interesting and valid results, a relatively
satisfactory result obtained by S.S.S.Wang(\cite{4}) in 1980 is that
JC holds for all polynomial mappings of degree 2 in all dimensions.
The most powerful and surprising result is the reduction to degree
3, due to H.Bass, E.Connell and D.Wright( \cite{2}) in 1982 and
A.Yagzhev(\cite{5}) in 1980, which asserts that JC is true if JC
holds for all polynomial mappings of degree 3(what is more, if JC
holds for all cubic homogeneous polynomial mapping!). In the same spirit of the above degree reduction method, another efficient way to tackle JC is the Druzkowski's Theorem(\cite{6}): JC is true if it is true for all Druzkowski mappings (in all dimension $\geq 2$). Another interesting result is due to Gorni-Zampieri (\cite{7}), who proved in 1997 that there exists Gorni-Zampieri pairing between the cubic homogeneous polynomial mappings and the Druzkowski mappings. One more interesting result is due to J.Gwozdziewicz (\cite{8}), who
proved in 1993 that any two-dimensional polynomial mapping injective
on one line is injective on the whole plane hence is an
automorphism, another beautiful result due to S.Abhyankar and
T.T.Moh in 1975(\cite{9}), which claims that an injective polynomial
mapping is indeed an automorphism (among the vast proofs given to
this result are M.Kang's(\cite{10}) and W.Rudin's(\cite{11}). In 1982, H.Bass E.Connell and D.Wright showed the degree of the inverse of polynomial automorphisms are no more than ${\operatorname{(degF)}^{n-1}}$ where ${\operatorname{degF}}=\max_{i}{\operatorname{degF_i}}$, that is equivalent ${\operatorname{degF}^{-1}}\leq {\operatorname{(degF)}^{n-1}}$(see [2]).

In this paper, we will prove that the homogeneous polynomials which satisfy the Jacobian condition are injective on the lines that pass through the origin. Therefore, the Jacobian Conjecture is true for $n=2$ in this case. In the second conclusion, we will show that $F$ and $G'$ are Gorni-Zampieri pairings, therefore, $G'$ is an isomorphism or $|JG'|=1$ if and only if $F$ has the same properties. Finally we will give a bound of ${\operatorname{degF}^{-1}}$ according to the rank of the matrix which consists of the coefficients of linear part of $F$.  \

\indent Recall that $F$ is a cubic homogeneous map if $F=X+H$ with
$X$ the identity (written as a column vector) and each component of
$H$ being either zero or cubic homogeneous. A cubic homogeneous
mapping $F=X+H$ is a {\it\bf Druzkowski (or cubic linear) mapping}
if each component of $H$ is either zero or a third power of a linear
form. Each Druzkowski mapping $F$ is associated to a scalar matrix
$A$ such that $F=X+(AX)^{*3}$, where $(AX)^{*3}$ is the {\it\bf
Druzkowski symbol} for the vector $(A^1X)^3,\cdots, (A^nX)^3)$ with
$A^i$ the $i-$th row of $A$. Clearly, a Druzkowski mapping is
uniquely determined by this matrix $A$. Apparently, the notion of a
Druzkowski mapping can be easily generalized. Namely, for any
positive integer $d\geq 2$, we call $F=X+H$ to be {\it\bf a
generalized Druzkowski mapping} of degree $d$ if each component of
$H$ is either zero or a $d$-th power of a
linear form.  A fact about a generalized Druzkowski
mapping $F$ is that if $|JF|\in \bf{C}^*$, then $JH$ is nilpotent
and $|A|=0$. Therefore, ${\operatorname{rankA}}\leq {n-1}$ . In fact, $JH$ is nilpotent that is equivalent to $|JF|\in \bf{C}^*$ in such case. Now, we can state the results of this paper:\\
\indent

{\bf Theorem 1 } {\it Let $F=X+H$, where H are homogeneous polynomials for ${\operatorname{degH=d}}$ and {\bf K} denotes to a character zero field. If ${\operatorname{detJF}}=1$. Then $F$ is injective on the lines that pass through the origin.}\\
\indent

{\bf Proof } Since $H=d^{-1}\Sigma_{j=1}^{n}x_{j}H_{x_{j}}$, therefore, we have $F=(I+d^{-1}JH)X$. For any $a\in {\bf K}^n, b=\lambda a$ for $\lambda \in {\bf K}$, suppose $F(a)=F(b)$, then we have\
$$(I+d^{-1}JH_{a})a=(I+d^{-1}JH_{b})b$$
where $JH_{a}, JH_{b}$ are got from $JH$ by replacing X by a, b respectively. So we have\
$$(I+d^{-1}JH_{a})a=(\lambda I+d^{-1}\lambda^dJH_{a})a$$\ this is equivalent
$$[(\lambda-1)I+d^{-1}(\lambda^d-1)JH_{a}]a=0$$ we have
$$(\lambda-1)[I+d^{-1}(\lambda^{d-1}+\cdots +1)JH_a]a=0$$ since $det(I+d^{-1}(\lambda^{d-1}+\cdots +1)JH_a)=1$, therefore, we have $\lambda=1$ or $a=0$. In any cases, we have $a=b$. This completes the proof.\\
\indent
In the proof of Theorem 1, we can get $X=(I+d^{-1}JH)^{-1}F$. Therefore, the ideal generated by $X$ is the same as the ideal generated by $F$, that is $(x_1, x_2,\cdots, x_n)=(F_1, F_2,\cdots, F_n)$.\\
\indent

{\bf Lemma2 } {\it Let {\bf K} be a field and assume $F:\bf{K}^n\rightarrow \bf{K}^n$ is an invertible polynomial map. If ${\operatorname{kerJ(F-x)}}\cap K^n$ has dimension $n-r$ as a K-space, then there exists a $T\in GL_n(K)$ such that for $G:=T^{-1}F(Tx)$, we have $kerJ(G-x)\cap \bf{K}^n={0}^r\times \bf{K}^{n-r}$. Furthermore, $G$ is invertible and the degree of the inverse is the same as that of $F$.}\\

\indent
{\bf Proof }Take $T\in GL_n(K)$ such that the last $n-r$ columns of $T$ are a basis of the K-space ${\operatorname{kerJ(F-x)}}\cap K^n$, and set $G:=T^{-1}F(Tx)$. Then
$$G-x=T^{-1}F(Tx)-T^{-1}Tx=T^{-1}(F-x)\mid_{x=Tx}$$
and by the chain rule, $J(G-x)=T^{-1}\cdot(J(F-x))\mid_{x=Tx}\cdot T$. Hence\\
\begin{eqnarray*}
   {\operatorname{kerJ(G-x)}}\cap K^n&=&T^{-1}({\operatorname{kerJ(F-x)}}\mid_{x=Tx})\cap K^n  \\
   &=&T^{-1}({\operatorname{kerJ(F-x)}})\mid_{x=Tx}\cap T^{-1}(K^n)\mid_{x=Tx}  \\
   &=&T^{-1}({\operatorname{kerJ(F-x)}}\cap K^n)\mid_{x=Tx}  \\
   &=&T^{-1}(KTe_{r+1})+KTe_{r+2}+\cdots+KTe_n)\mid_{x=Tx}  \\
   &=&Ke_{r+1}+Ke_{r+2}+\cdots+Ke_n  \\
   &=&\{0\}^r\times K^{n-r}
\end{eqnarray*}
If \~{F} is the inverse of $F$, then \~{G}$:=T^{-1}$\~{F}(Tx) is the inverse of $G$ and \~{G} have the same degree as \~{F}, which completes the proof.\\
\indent

{\bf Theorem3 } {\it Let $K$ be a field with $charK=0$ and assume $F=x+H:K^n\rightarrow K^n$ is an invertible polynomial map of degree $d$. If ${\operatorname{kerJH}}\cap \bf{K}^n$ has dimension $n-r$ as a K-space, then the inverse polynomial map of F has degree at most $d^r$.}\\

\indent
{\bf Proof }From lemma2, it follows that by replacing $F$ by $T^{-1}F(Tx)$ for a suitable $T\in GL_n(K)$, we may assume that ${\operatorname{kerJH}}\cap \bf{K}^n=\{0\}^r\times \bf{K}^{n-r}$. This means that the last $r$ columns of $JH$ are zero. Hence $H_i\in K[x_1, x_2,\cdots,x_r]$ for all $i$. Let $x-G$ be the inverse of $F$ at this stage. Then $G(F)=H$, since $x+H-G(x+H)=x$.

Let $L=K(x_1,x_2,\cdots,x_r)$. Since $F_i, H_i\in L$ for all $i\leq r$ and $F_{r+1}, F_{r+2},\cdots,F_n$ are algebraically independent over $L$, we obtain from $G_i(F)=H_i$ that $G_i\in K[x_1,x_2,\cdots,x_r]$ for all $i\leq r$. Hence $(x_1-G_1,x_2-G_2,\cdots,x_r-G_r)$ is the inverse polynomial map of $(F_1,F_2,\cdots,F_r)$. Therefore, $deg(x_1-G_1,x_2-G_2,\cdots,x_r-G_r)\leq d^{r-1}$ on account of [3,Prop2.3.1].

Using that F is the inverse of $x-G$, we obtain by substituting $x=x-G$ in $G_i(F)=H_i$ that $G_i=H_i(x-G)$ for all $i$. Since $degH_i\leq d, H_i\in K[x_1,x_2,\cdots,x_r]$ and ${\operatorname{deg(x_1-G_1,x_2-G_2,\cdots,x_r-G_r)}}\leq d^{r-1}$, we get ${\operatorname{degG_i}}\leq d^r$ for all $i$. This completes the proof.
\indent

{\bf Corollary4 } {\it Let {\bf K} be a field with $charK=0$ and assume $F:{\bf K}^n\rightarrow {\bf K}^n$ is an invertible polynomial map of the form $F=x+H$, where $H$ is power linear of degree d. If ${\operatorname{rankJH}}=r$, then the inverse polynomial map of $F$ has degree at most $d^r$.}\\

\indent
{\bf Proof }Since $H=(Ax)^{*d}$. Then we have $JH=d\cdot {\operatorname{diag((Ax)^{*(d-1)})}}\cdot A$, whence ${\operatorname{kerJH}}={\operatorname{kerA}}$ ${\operatorname{rankA}}={\operatorname{rankJH}}=r$. Since $A$ is a matrix over $K$, the dimension of ${\operatorname{kerA}}\cap {{\bf K}}^n$ as a K-space is $n-r$. Hence the desired result follows from theorem 3.\\\

 Acknowledgment: The author is very grateful to Professor Yuehui Zhang who introduced the Conjecture and gave great help when the author studied the problem and Professor Guoping Tang who read the paper carefully and gave some good advice. The author also wishes to thank the anonymous referee, whose valuable comments and remarks helped to improve the paper.

\end{document}